\documentclass[12pt]{article}

\usepackage{anysize}
\marginsize{28mm}{28mm}{25mm}{25mm}
\usepackage{setspace}
\singlespace

\usepackage[dvips]{graphicx}
\usepackage[fleqn]{amsmath}

\usepackage{amsmath,amssymb}
\usepackage{epsfig,subfigure}
\usepackage{algorithmic}
\usepackage{stmaryrd}
\usepackage{txfonts}
\usepackage{multirow}
\usepackage{url}
\usepackage[subnum]{cases}

\usepackage{algorithmic}
\usepackage{algorithm}

\usepackage{booktabs}

\usepackage{pifont}

\input{def1.set}

\input{psfig.sty}

\makeatletter
  \newcommand{\figcaption}[1]{\def\@captype{figure}\caption{#1}}
  \newcommand{\tblcaption}[1]{\def\@captype{table}\caption{#1}}
\makeatother

\begin{document}

\title{Tensor Decompositions: A New Concept in Brain Data Analysis?}

\author{Andrzej CICHOCKI
\thanks{a.cichocki@riken.jp, Laboratory for Advanced Brain Signal Processing, RIKEN BSI Japan,
Please also see  \bf {Control Measurement, and System Integration (SICE),
special issue; Measurement of Brain Functions and Bio-Signals}, 7, 507-517, (2011).}}

\maketitle

{\bf Abstract}

Matrix factorizations  and their extensions to tensor factorizations and decompositions have become  prominent techniques for linear and
multilinear blind source separation (BSS), especially multiway Independent Component Analysis (ICA), Nonnegative Matrix  and Tensor Factorization (NMF/NTF), Smooth Component Analysis (SmoCA) and Sparse Component Analysis (SCA). Moreover, tensor decompositions have many other potential applications beyond multilinear BSS, especially feature extraction, classification, dimensionality reduction  and multiway clustering.
In this paper, we  briefly  overview new and emerging models and approaches for tensor decompositions in applications to  group  and linked multiway BSS/ICA, feature extraction, classification  and Multiway Partial Least Squares (MPLS) regression problems.


{\bf Keywords:} Multilinear BSS, linked multiway BSS/ICA, tensor factorizations and decompositions, constrained Tucker and CP models, Penalized Tensor Decompositions (PTD),  feature extraction, classification, multiway PLS and CCA.

\section{Introduction}
\label{sec:introduction}
 Although the basic models for tensor (i.e.,  multiway array) decompositions and factorizations  such as Tucker and (Canonical-decomposition/Parafac) CP models were proposed long time ago,  they are recently emerging as promising tools for exploratory analysis of multidimensional data in diverse applications, especially, in multiway blind source separation (BSS), feature extraction, classification, prediction, and multiway clustering \cite{NMF-book,Kolda08,Lath-BCM,Hitchcock1927,Tucker66}. The recent advances in neuroimage technologies (e.g., high density array EEG/MEG, fMRI, NIRS) have generated massive amounts of brain data exhibiting high dimensionality, multiple modality and multiple couplings, functional connectivity. By virtue of their multiway nature, tensors provide a powerful tools for analysis and fusion of such massive data together with a mathematical backbone for the discovery of underlying hidden  complex data structures \cite{NMF-book,Cichocki-GCA,CiAm02}.
%
%
Constrained Matrix Factorizations, called sometimes penalized matrix decompositions (e.g., ICA, SCA, NMF, SVD/PCA, CCA) and their extensions or penalized tensor decompositions - multidimensional constrained models  (Tucker,  CP, NTF, NTD models \cite{Cichocki-GCA,NMF-book,CiAm02,Phan-HALS-IEICE}), with some constraints such as statistical independence, decorrelation, orthogonality, sparseness, nonnegativity, and smoothness - have been recently proposed as meaningful and efficient representations of signals, images and in general natural multidimensional data \cite{NMF-book}.
 From signal processing and  data analysis point of view,  tensor decompositions
 are very attractive because they take  into account spatial, temporal and spectral information and provide links among various data and extracted  factors or hidden
(latent variables)  components with  physical or physiological meaning and
interpretations  \cite{NMF-book,Phan-HALS-IEICE,Morup,Acar-Morup11,Miwakeichi}.
In fact, tensor decompositions are emerging techniques for data fusion, dimensionality reduction, pattern recognition, object detection, classification, multiway clustering, sparse representation and coding, and multilinear blind source separation (MBSS) \cite{NMF-book,Kolda08,CichZd_ICASSP07,Phan-Class10,Lathauwer_HOOI}.

{\bf  Basic notations.} Tensors (i.e., multiway arrays) are denoted by underlined capital boldface letters, e.g., $\underline \bY \in \Real^{I_{1} \times I_{2} \times \cdots \times I_{N}}$.
The order  of a tensor is the number of modes, also known as ways or dimensions (e.g., space, time, frequency, subjects, trials, classes, groups, conditions).
In contrast, matrices (two-way tensors) are denoted by boldface capital letters, e.g., $\bY$;
 vectors (one-way tensors) are denoted by boldface lowercase letters, e.g., columns of the matrix $\bU$ by $\bu_j$
 and scalars are denoted by lowercase letters, e.g.,
$u_{ij}$.

The mode-$n$  product $\underline \bY=\underline \bG \times_n \bU$ of a tensor
$\underline \bG \in \Real^{J_{1} \times J_{2} \times \cdots \times J_{N}}$ and
a matrix $\bU \in \Real^{I \times J_n}$ is  a tensor
$\underline \bY \in \Real^{J_1 \times \cdots \times J_{n-1} \times I \times J_{n+1} \times \cdots \times J_N}$, with elements $y_{j_1,j_2,\ldots,j_{n-1},i_n,j_{n+1},\ldots, j_{N}} =\sum_{j_n=1}^{J_n} g_{j_1,j_2,\ldots,J_N} \; u_{i_n,j_n}$.
Unfolding (matricization, flattening) of a tensor  $\underline \bY \in \Real^{I_{1} \times I_{2} \times \cdots \times I_{N}}$ in $n$-mode is denoted as $\bY_{(n)} \in \Real^{I_{n} \times (I_{1} \cdots I_{n-1} I_{n+1}, \cdots I_N)}$, which consists of arranging all possible $n$-mode tubes in  as columns of a matrix. \cite{Kolda08}.
Throughout this paper, standard notations and basic tensor operations are used \cite{NMF-book}.

\section{Basic Models for  Multilinear BSS/ICA}

Most of the linear blind source separation (BSS) models can be represented as  constrained matrix factorization problems, with suitable constraints imposed on factor matrices (or their columns -referred to as components)
\be
 \bY = \bA \bB^T + \bE =\sum_{j=1}^J \ba_j \bb^T_j +\bE =\sum_{j=1}^J \ba_j \circ \bb_j +\bE,
\label{NMF1}
 \ee
where
$\circ$ denotes outer product\footnote{The outer product of two vectors
$\ba \in \Real^I, \; \bb \in \Real^T$
builds up a rank-one matrix
$\bY=\ba \circ \bb = \ba \bb^T \in \Real^{I \times T}$
and the outer product of  three vectors: $\ba \in \Real^I, \; \bb \in \Real^T, \; \bc \in \Real^Q$ builds up a third-order rank-one tensor:
$\underline \bY = \ba \circ \bb \circ \bc \in \Real^{I \times T  \times Q}$,
with entries defined as $y_{itq} =a_i \; b_t \; c_q$.},
$\bY = [y_{it}] \in \Real^{I \times T}$ is a known  data matrix (representing observations or measurements),
$\bE = [e_{it}]  \in \Real^{I \times
T}$ represents errors or noise, $\bA=[a_{ij}]=[\ba_1,\ba_2,$ $\ldots,\ba_J]$ $ \in \Real^{I \times J}$ is an unknown  basis (mixing) matrix, with basis vectors $\ba_{j} \in  \Real^I $ and the $J$ columns of a matrix   $\bB = [ \bb_1, \bb_2, \ldots, \bb_J] \in \Real^{T \times J}$
represent unknown components. latent variables or sources  $\bb_j$.

{\bf Remark}:We notice that we have
symmetry in the factorization: For Eq (\ref{NMF1}) we could just as easily write $\bY^T \cong \bB \bA^T$,  so the meaning of "sources" and "mixture" are often somewhat arbitrary.

 Our primary objective in the BSS is to estimate uniquely (neglecting unavoidable scaling and permutation ambiguities) factor matrices $\bA$ and $\bB$, subject to   various  specific constraints imposed on the vectors $\bb_j$ and/or $\ba_j$,  such as mutual statistical independence (ICA), sparseness (SCA), smoothness (SmoCA),  nonnegativity (NMF) or  orthogonality (PCA/SVD), uncorrelation, etc.

In some applications the data matrix  $\bY$ is factorized into three or more factors \cite{NMF-book}.
In the special case, of a singular value decomposition (SVD) of the data matrix $\bY \in \Real^{I \times T}$, we have the following factorization:
\be
\bY= \bA \bD \bB^T = \bD \times_1 \bA \times_2 \bB = \sum_j d_{jj} \ba_j \bb_j^T,
\label{SVD1}
\ee
where $\bA \in \Real^{I \times I}$ and $\bB\in \Real^{T \times T}$ are orthogonal matrices and $\bD$ is a diagonal matrix  containing only nonnegative singular values. The SVD and its generalizations play key roles in signal processing and data analysis \cite{Lathauwer_HOOI}.

%

Often  multiple subject, multiple task data sets can be represented by a set  of data matrices $\bY_n$ and it is necessary to perform simultaneous constrained matrix
factorizations\footnote{This form of factorizations are typical with EEG/MEG related data for multi-subjects, multi-tasks, while the factorization for the transposed of $\bY$ is typical
for fMRI data \cite{GICA-rev,Guo-GICA,Langers-GICA}.}:
\be
\bY_n \cong \bA_n \bB_n^T \qquad (n=1,2,\ldots,N)
\label{GBSS1}
\ee
or in preprocessed form with a dimensionally reduction
\be
 \tilde\bY_n \cong \bQ_n \tilde \bA_n \bP_n \tilde \bB_n^T, (n=1,2,\ldots,N),
\label{GBSSr}
\ee
subject to various additional constraints (e.g., $\bB_n=\bB, \; \forall n $ and their columns are mutually independent and/or sparse). This problem is related to various models of group ICA, with suitable pre-processing,  dimensionality reduction and post-processing procedures \cite{GICA-rev,Guo-GICA,Langers-GICA}.
In this paper, we introduce  the group multiway BSS concept, which is more general and flexible
 than the  group ICA, since  various constraints can be imposed on factor matrices in different modes (i.e., not only mutual independence but also nonnegativity, sparseness, smoothness or orthogonality).
%
There is neither a theoretical nor an experimental basis that statistical independence (ICA) is the unique right concept to extract brain sources or isolate brain networks \cite{Cichocki-GCA}.
In real world scenarios, latent (hidden) components (e.g., brain sources) have various
complex properties and features. In other words, true unknown involved components are seldom
all statistically independent. Therefore, if we apply only one single criterion like ICA, we may fail to extract all desired  components with physical interpretation.
We need rather to apply an ensemble  of strategies by employing several suitably chosen criteria and associated learning algorithms to extract all desired components  \cite{NMF-book,Cichocki-GCA}. For these reasons, we have  developed multiway BSS methods, which are based not only on the statistical independence of components, but rather  exploits multiple criteria or diversities of factor matrices in various modes,  in order to extract physiologically meaningful  components, with specific  features and statistical properties
\cite{NMF-book,Cichocki-GCA}. By diversity, we mean  different characteristics, features or morphology of source signals or hidden latent variables \cite{CiAm02}.
%
Since multi-array data can be always interpreted in many different ways, some {\it a priori}
knowledge is needed to determine, which diversities, characteristics, features or properties represent  true latent (hidden) components with physical meaning.

It should be noted, although standard 2D BSS (constrained matrix factorizations) approaches, such as ICA, NMF, SCA, SmoCA, PCA/SVD, and their variants, are invaluable tools for feature extraction and selection, dimensionality reduction, noise reduction, and data mining, they have only two modes or 2-way representations (typically, space and time), and their use is therefore limited. In many neuroscience applications the data structures often contain higher-order ways (modes) such as subjects,  groups, trials, classes and conditions, together with the intrinsic dimensions of space, time, and frequency.
For example, studies in neuroscience often involve multiple subjects (people or animals) and trials leading to experimental data structures conveniently represented by multiway arrays or blocks of multiway data.

%

\begin{figure}[t!]
\begin{center}
(a)
\includegraphics[width=10.2 cm,height=4.2cm]{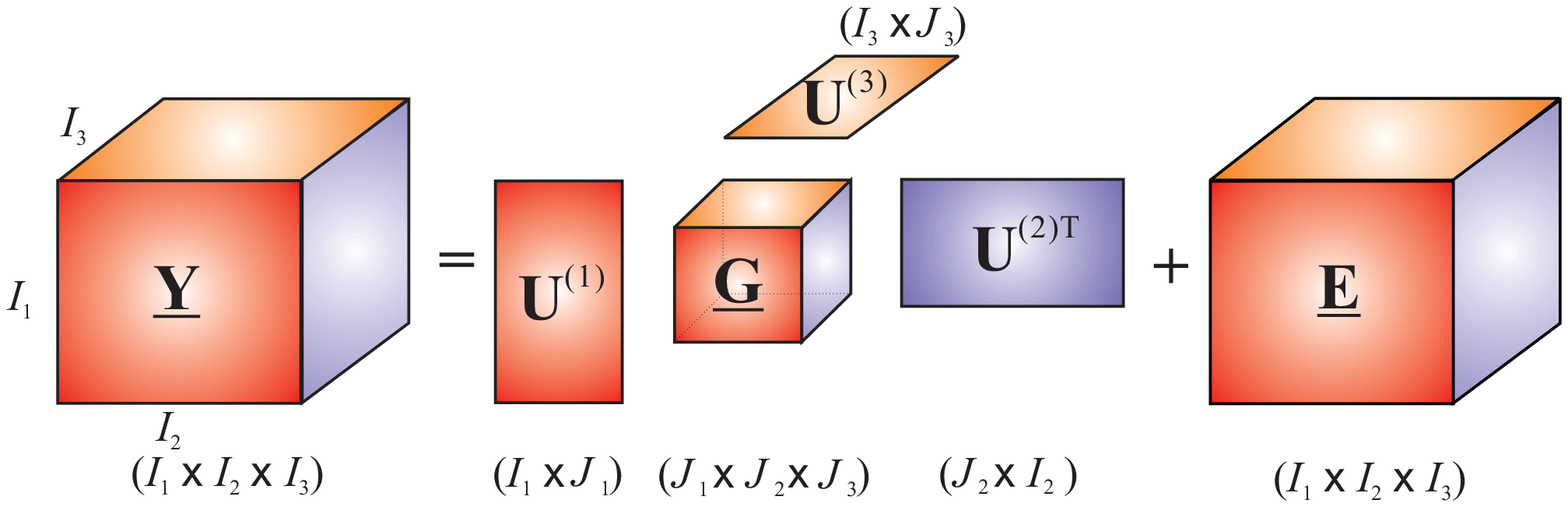}\\
\vspace{0.7cm}
(b)
\includegraphics[width=10.4 cm,height=3.1cm]{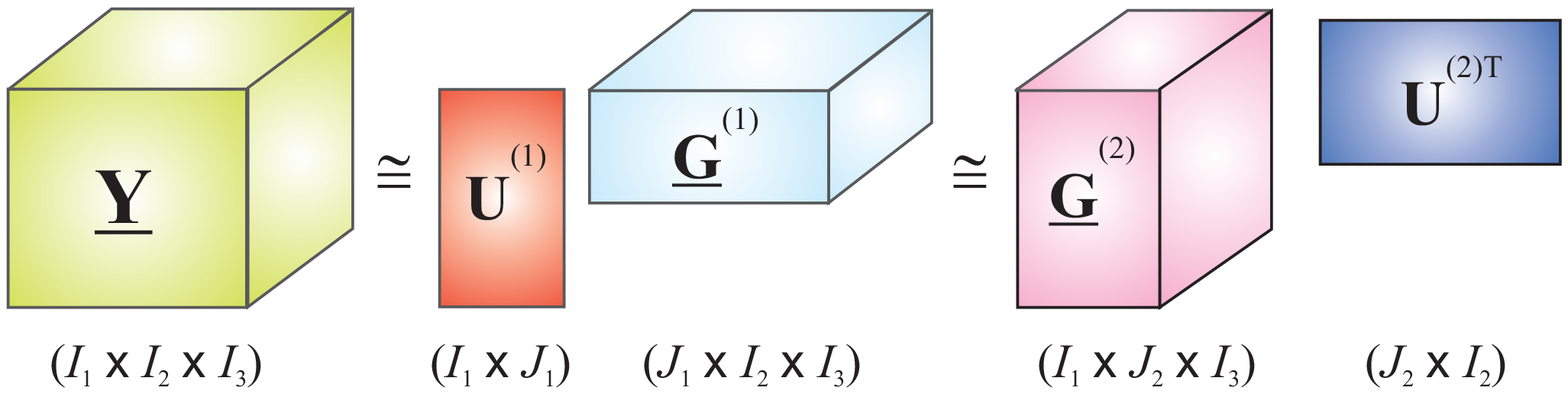}
\end{center}
\caption{Illustration of a 3-way tensor decomposition using a constrained Tucker-3 model; Objective is to estimate  factor matrices $\bU^{(n)}=[\bu_1^{(n)},\bu_2^{(n)},\ldots,\bu_{J_n}^{(n)}] \in \Real^{I_n \times J_n}$  (with desired diversities or statistical properties) and a possibly sparse core tensor $\underline \bG \in \Real^{J_1 \times J_2 \times J_3}$, typically with $J_n << I_n$ $(n=,1,2,3)$. (b) Preselected factor matrices in Tucker-3 can be absorbed by a core tensor $\underline \bG$, and this leads to Tucker-1 models: $\underline \bY = \underline \bG^{(1)} \times_1 \bU^{(1)}+\underline \bE^{(1)}= \underline \bG^{(2)} \times_2 \bU^{(2)}+\underline \bE^{(2)}=\underline \bG^{(3)} \times_3 \bU^{(3)}+\underline \bE^{(3)}$. Instead of applying the standard Alternating Least Squares (ALS) algorithms for the Tucker-3 model, we can apply unfolding of data tensor according to the Tucker-1 models and then perform constrained factorizations of the unfolded matrices (multiway BSS)  by  imposing desired constraints (independence, sparseness, nonnegativity, smoothness or uncorrelation, etc.)}.
\label{fig:Tucker3}
\vspace{-4mm}
\end{figure}

The simple linear BSS models (\ref{NMF1})-(\ref{GBSS1}) can be naturally extended for multidimensional data to the multiway BSS models using constrained tensor decompositions.
In this paper, we consider a general and flexible approach  based on  the Tucker decomposition model - called also Tucker-N model (see Fig. \ref{fig:Tucker3} (a)) \cite{Tucker66,NMF-book,Kolda08}:
\be
\underline{\bY}
&=&
 \sum\limits_{j_1 = 1}^{J_1} {
 \sum\limits_{j_2 = 1}^{J_2} {\cdots
 \sum\limits_{j_N = 1}^{J_N}{g_{j_1 j_2 \cdots j_N} \, \left(\; \bu^{(1)}_{j_1} \circ \bu^{(2)}_{j_2} \circ \cdots \circ \bu^{(N)}_{j_N} \right) }}} + \underline{\bE} \notag\\
&=& \underline{\bG} \times_1 \bU^{(1)}\times_2 \bU^{(2)}\cdots \times_N \bU^{(N)} + \underline{\bE} \notag\\
&=& \underline{\bG} \times \{ \bU \} + \underline{\bE} = { \underline {\widehat \bY}} + \underline{\bE},
\label{TuckerN}
\ee
where $\underline \bY \in \Real^{I_1 \times I_2 \cdots \times I_N}$ is the given data tensor, $\underline \bG \in \Real^{J_1 \times J_2 \cdots \times J_N}$ is a core tensor of reduced dimension, $\bU^{(n)}=[\bu^{(n)}_{1}, \bu^{(n)}_{2},\ldots,\bu^{(n)}_{J_n}]\in \Real^{I_n \times J_n}$ ($n = 1, 2,\ldots, N$) are factors (component matrices) representing components, latent variables, common factors or loadings, ${\underline{\widehat {\bY}}}$ is an approximation of the measurement $\underline{\bY}$, and $\underline{\bE}$ denotes the ap\-proxi\-ma\-tion error or noise depending on the context \cite{NMF-book,Tucker66,Kolda08}.
The objective is to estimate  factor matrices: $\bU^{(n)}$, with components (vectors) $\bu^{(n)}_{j_n}$, $(n=1,2,\ldots, N, \;\;j_n=1,2,\ldots,J_N)$ and the core tensor $\underline \bG \in \Real^{J_1 \times J_2 \cdots \times J_N}$  assuming  that the number of factors in each mode  $J_n$ are known or can be estimated \cite{NMF-book}.

If the factor matrices and a core tensor are orthogonal the Tucker model can be considered as extension of the SVD model (\ref{SVD1}), known as the High Order SVD (HOSVD)  \cite{Lathauwer_HOOI}. While
the optimal approximation of a matrix can be obtained by truncation
of its SVD, the optimal tensor approximation  (with a minimal norm of the error) cannot in general
be obtained by truncation of the Tucker decomposition. However, it was shown that the
truncation of the particular constrained version  usually yields a
quite  good approximation \cite{Lathauwer_HOOI}. 
 In general, in the Tucker model the orthogonality
constraints are usually not  imposed. Instead, we consider alternative constraints, such as sparseness, nonnegativity, smoothness  or statistical mutual independence.
It should be noted that the Tucker model is not unique in general, even if we impose some weak constraints. However, this model is unique if we impose suitable constraints (e.g., sparseness or independence).

This leads to a concept of group multiway BSS. There are two possible interpretations or concepts of employing Tucker decomposition as multiway BSS. In the first concept the columns of factor matrices $\bU^{(n)}$ represent desired components or latent variables and the core tensor represent some "mixing process" or more precisely the core tensor shows links among components in different modes, while data tensor $\underline \bY$  represents collection of 1-D or 2-D mixing signals. In the second concepts, the core tensor represents desired but unknown (hidden)  $N$-dimensional signal (e.g., 3D MRI image or 4D video) and factor matrices represent various transformations, e.g., time frequency transformation or wavelets dictionaries (mixing  or filtering processes), while a data tensor $\underline \bY$ represents observed $N$-dimensional signal which is distorted, transformed compressed or mixed depending on applications. In this paper, we will consider only the first interpretation.

The Tucker-$N$ model (\ref{TuckerN}) can be represented by $n$ approximative matrix factorizations with three factors:
\be
\bY_{(n)} \cong \bU^{(n)} \; \bG_{(n)} \; \bZ^{(n)}, \qquad (n=1,2, \ldots,N),
\ee
where $\bZ^{(n)}=\left[\bU^{(N)} \otimes \cdots \otimes \bU^{(n+1)}
  \otimes \bU^{(n-1)}  \cdots \otimes \bU^{(1)}\right]^{T} $.

Moreover, the Tucker-$N$ model (\ref{TuckerN}) can be compressed to $N$ Tucker-1 models, with only one factor  matrix $\bU^{(n)}$, in each mode $n$ (see Fig \ref{fig:Tucker3} (b)):
\be
\underline \bY \cong \underline \bG^{(n)} \times_n \bU^{(n)} \;\;\mbox{\small or in matrix form} \;\; \bY_{(n)} \cong  \bU^{(n)} \bG_{(n)}^{(n)},
\label{Tucker1}
\ee
where $\underline \bG^{(n)}= \underline \bG \times_1 \bU^{(1)}\times_2  \cdots \times_{n-1} \bU^{(n-1)} \times_{n+1} \bU^{(n+1)}  \cdots \times_N \bU^{(N)}$, $\;(n=1,2,\ldots,N)$.
   Note  that the above models correspond
 to group BSS/ICA models (\ref{GBSS1}), with $\bY_n=\bY_{(n)}^T$, $\bA_n= [\bG^{(n)}_{(n)}]^T$ and $\bB_n=\bU^{(n)}$, $\;\forall n$, under the assumption that we impose desired constraints for factor matrices $\bU^{(n)}$.

 Most of the existing approaches exploit only the CP model and impose only statistical independence constraints This leads to tensor probabilistic ICA \cite{Beckmann2005} or ICA-CPA  models \cite{DeVos2008,LinkedICA}.  
 However, our approaches are quite different, because we use the Tucker models and we are not restricting only to ICA assumptions but exploit multiple criteria and allows for diversities of components. The advantage of the Tucker model over the CP model is that the  Tucker model is more general and the  number of components in each mode  can be different and furthermore the components are linked via a core tensor, and hence  allows us to model more complex hidden  data structures. Note that the Tucker decomposition can be simplified to the CP model in
the special case, where the "hyper-cube"   core tensor (with $J=J_1=J_2 = \cdots =J_N$) has nonzero elements only on the super-diagonal.
The CP model has unique factorization without any constraints under
some mild conditions.

For the Tucker and the  CP models there exist many efficient algorithms
\cite{NMF-book,Kolda08}. Most of them are based on the ALS (Alternating
Least Squares) and HALS (Hierarchical ALS) \cite{NMF-book,Phan-HALS-IEICE,Kopriva-Jeric-Ci,Acar-Morup11}
and CUR (Column-Row) decompositions \cite{Caiafa-Cichocki-CUR}.
However, description of these algorithms is out of the scope of this paper.

 We can implement Multilinear BSS algorithms in several ways.
First of all,  we can minimize a  global cost function,  with  suitable penalty and/or regularization terms to estimate desired components (see Eq. (\ref{TuckerN})):
\be
D_F\left( \underline{\bY}\| \underline{\bG}, \{\bU\}\right) = \| \underline{\bY} - \underline {\bG} \times \{\bU\} \|_F^2 
+ \sum_n \alpha_n C_n(\bU^{(n)}),
\label{cost-Tucker}
\ee
where $\alpha_n \geq0$ are penalty coefficients and $C_n(\bU^{(n)})$ are penalty terms, which are added to achieve specific characterstic of the components. For example, if we need to impose mutual independence constraints the penalty terms can take the following form $C_n(\bU^{(n)})= \sum_{n=1}^{N} \sum_{j\neq p} \bu^{(n)T}_{j} \psi_p (\bu^{(n)}_{n})$, where $\psi_n(u)$ are suitable nonlinear functions.
In principle, this method  referred as penalized tensor decompositions, allows to find the  factor matrices
$\bU^{(n)}$, with unique components $\bu_{j_n}^{(n)}$  and an associated core tensor
but the method involves heavy computations and it is time consuming.

Another approach is to apply standard Tucker decompositions method,  without applying any desired constraints using  ALS \cite{NMF-book},  CUR \cite{Caiafa-Cichocki-CUR} or HOOI/HOSVD  \cite{Lathauwer_HOOI,MICA2005} algorithms and in the next step apply for each factor matrix standard constrained matrix factorization  algorithms (e.g., ICA, NMF or SCA) \cite{ICATucker,Zhou-MBSS}.
Assuming that each unconstrained factor $\bU^{(n)}$ in the model (\ref{TuckerN}) can be further factorized using standard BSS algorithms as $\bU^{(n)} \cong \tilde \bB_n \tilde \bA_n^T$,
 we can formulate the following decomposition model:
\be
{\underline \bY} & \cong& {\underline \bG} \,
 \times_1 \bU^{(1)} \, \times_2 \bU^{(2)} \cdots \, \times_N \bU^{(N)}  \\
  &\cong& \widetilde {\underline \bG} \,
 \times_1 \tilde \bB_{1} \, \times_2 \tilde \bB_{2} \cdots \, \times_N \tilde\bB_{N},
 \label{Tucker-ICA}
 \ee
 where $\widetilde {\underline \bG} =\left[{\underline \bG} \times_1 \tilde \bA_1^T \times_2 \tilde \bA_2^T \cdots \times_N \tilde \bA_N^T \right]$. It is worth to note, that in each mode, we can apply different criteria for matrix factorization.

 Alternatively, a simpler approach is to perform the unfolding the data tensor $\underline \bY$  for each mode $n$, according to the Tucker-1 decompositions (\ref{Tucker1})
  and to apply directly a suitable  constrained factorization of matrices (or a penalized matrix decomposition) \footnote{Unfortunately, this approach does not guarantee best fitness of the model to the data or minimum norm of errors $||\bE||_F^2$, but usually the solutions are close to optimal. Note that in practice, for large scale problems, we do not need to perform explicitly unfolding of the full data tensor. Instead, we may apply fast sampling of tubes (columns) of  data tensors \cite{Caiafa-Cichocki-CUR}.}
  \be
  \bY_{(n)} \cong \bB_n \bA_n^T \quad \mbox{or} \quad  \bY_{(n)}^T \cong \bA_n \bB_n^T, \quad (n=1,\ldots,N)
  \ee
  subject to desired  constraints,  by employing  standard efficient  BSS algorithms (e.g., NMF or ICA) \cite{NMF-book,Zhou-MBSS}.
  Since some  matrices $\bY_{(n)}$ may have large dimensions, we  usually need to apply
  some efficient methods for dimensionality reduction \cite{Zhou-MBSS}. Finally, the core tensor, which shows the links among components in different modes can be computed as
\be
 \widetilde{\underline \bG} = {\underline \bY} \,
 \times_1 \left[\bB_{1}\right]^{+} \, \times_2 [\bB_{2}]^+ \cdots \, \times_N [\bB_{N}]^+,
 \label{ProjG}
 \ee
 where $[\cdot]^+$ denotes Moore-Penrose pseudo-inverse of a matrix.

 Finally, we can apply the Block Tensor Decomposition (BTD) approach \cite{Lath-BCM}, with suitable constraints imposed on all or only in some preselected factors or components. In this particular  case, the simplest scenario is to employ a constrained Block Oriented Decomposition (BOD) model (combining or averaging $N$ Tucker-1 models) \cite{NMF-book}:
 \be
\underline \bY = \frac{1}{N} \sum_{n=1}^N \left(\underline \bG^{(n)} \times_n \bU^{(n)} \right)+ \underline {\tilde \bE}.
\ee
This model allows us to estimate all desired components in parallel or sequentially by minimizing the norm of the total errors ($||\underline {\tilde \bE}||_F^2$) subject to specific constraints.

\section{Dimensionality Reduction, Feature Extraction and Classification of Multiway Data}
\label{sec::classificaion2D}

\begin{figure}[t!]
 \centering
 (a) \vspace{0.2cm}\\
 \includegraphics[width=8.5 cm,height=7.8cm]{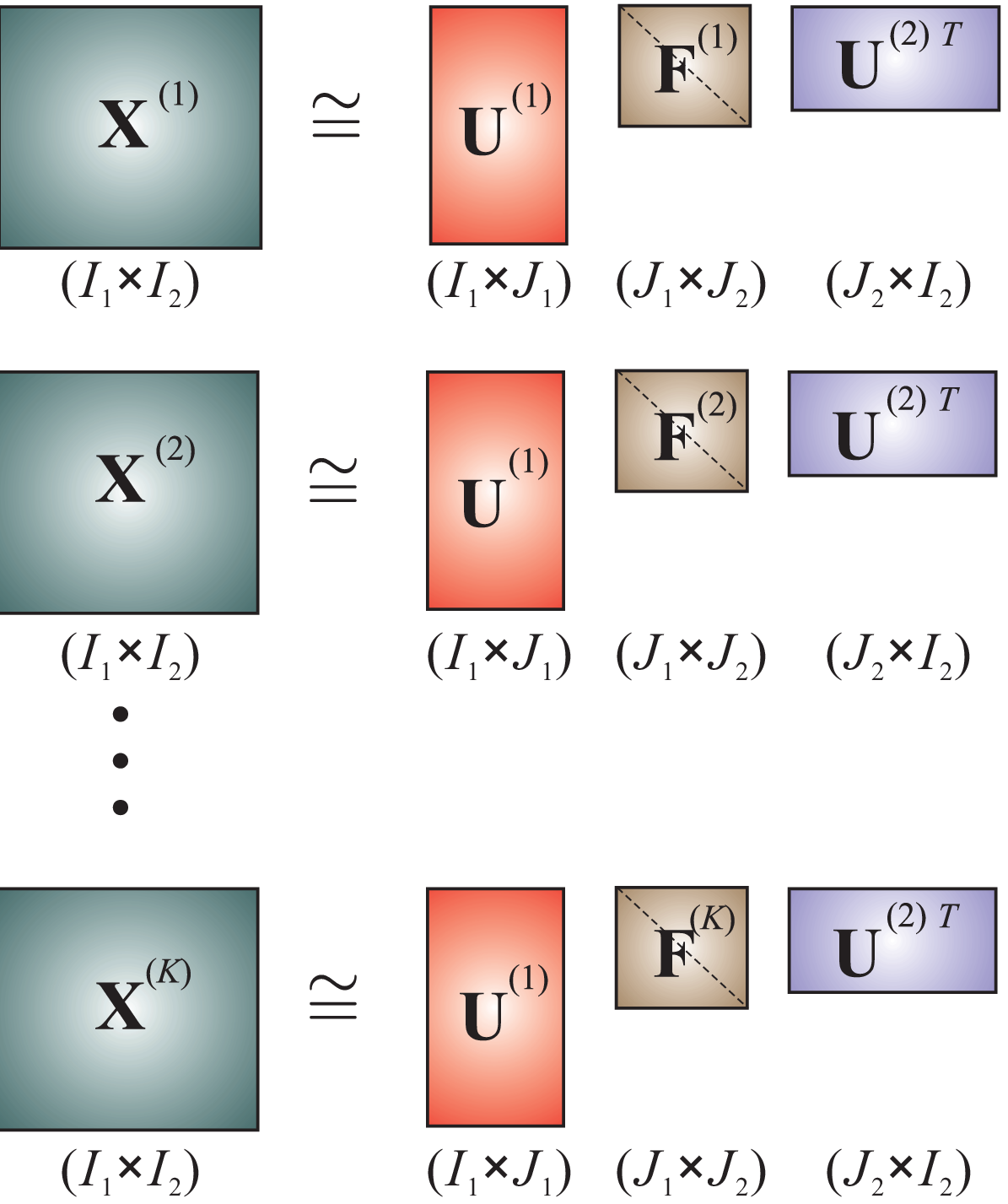} \vspace{0.7cm}\\
 (b)  \vspace{0.2cm} \\
 \includegraphics[width=10.5 cm,height=3.2cm]{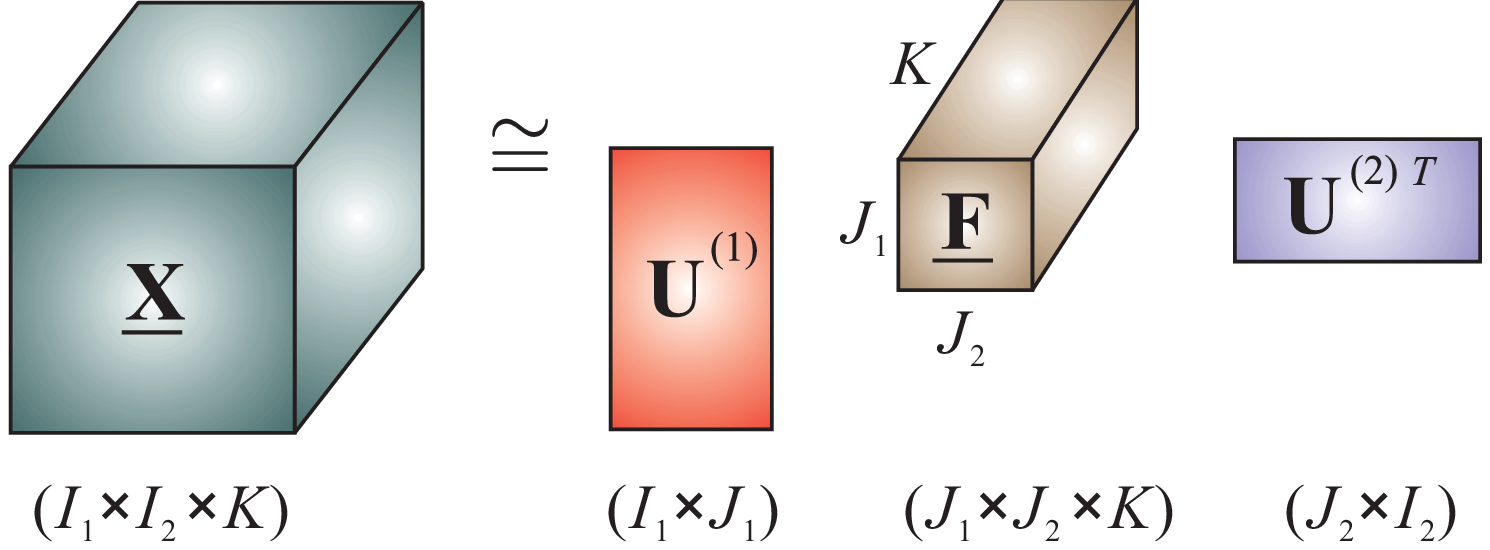}
 \caption{(a) Simultaneous matrix factorizations for dimensionality reduction and feature extraction. (b) This problem  is equivalent to a Tucker-2
 decomposition of a 3-way tensor into a constrained core tensor
 $\underline{\mathbf{F}} \in \Real^{J_1 \times J_2 \times K}$
 (representing features) and two basis factors
 $\mathbf{U}^{(1)} \in \Real^{I_1 \times J_1}$ and $\mathbf{U}^{(2)} \in \Real^{I_2 \times J_2}$.}
 \label{fig:FA1}
\end{figure}

Dimensionality reduction, feature extraction and selection are essential problems in the analysis of multidimensional datasets with large number of variables \cite{Phan-Class10}.
We shall first illustrate the basic concepts of dimensionality reduction
and feature extraction for  a
set of large-scale sample matrices.  Let us consider, that we have available
a set of $K$  matrices (2-D samples) ${\bX}^{(k)} \in
\Real^{I_1 \times I_2}, (k = 1,2, \ldots, K)$ that  represent multi subjects or multi-trials
2D data, which belong to $C$ different classes or categories (e.g., multiple mental task/state data or different mental diseases).  In order to perform dimensionality reduction and to extract
essential features for all the training samples, we apply simultaneous (approximative and constrained) matrix factorizations (see  Fig.~\ref{fig:FA1} (a)):
\be
\bX^{(k)} = \bU^{(1)} \, \bF^{(k)} \,
\bU^{(2)\, T} +  \bE^{(k)}, \qquad (k = 1,2,\ldots, K), \label{equ_2D_factorization}
\ee
where the two common factors (basis matrices) $\bU^{(1)} \in
\Real^{I_1 \times J_1}$ and $\bU^{(2)} \in \Real^{I_2 \times J_2}$, $J_n
\leq I_n,\;\forall n$ code (explain) each sample $\bX^{(k)}$ simultaneously along the
horizontal and vertical dimensions and the extracted features are
represented by matrices $\bF^{(k)} \in \Real^{J_1 \times J_2}$,
typically, with $J_1 << I_1$ and $ J_2<< I_2$. In special cases  $\bF^{(k)}$ are squared diagonal
matrices. This problem can be considered as a generalization of Joint Approximative Diagonalization (JAD) \cite{NMF-book,deLathauwer-JAD}.

The common method to solve such matrix factorizations problem is to
minimize a set of cost functions $||\bX^{(k)} - \bU^{(1)} \, \bF^{(k)}
\, \bU^{(2)\, T}||^2_F, \;\; \forall k$ sequentially or in parallel, with respect to all the factor
matrices.
We can solve the problem more efficiently by
concatenation or tensorization  of all matrices $\bX^{(k)}$ along the third dimension
to form an $I_1 \times I_2 \times K$ dimensional data tensor
$\underline{\bX} \in \Real^{I_1 \times I_2 \times K}$ and perform the standard Tucker decomposition (see Fig. \ref{fig:FA1} (b)).
%

In order to obtain meaningful components and a unique decomposition
 it is often convenient to impose additional constraints \cite{NMF-book,Phan-Class10}.
 In the special case when the feature matrices $\bF_k$ are diagonal the Tucker-2 model
 is reduced to special form of the unique CP model \cite{NMF-book}.

This approach can be naturally and quite straightforwardly  extended to multidimensional data.
Assume, that we have available a set of multidimensional tensors $\underline{\bX}^{(k)} \in \Real^{I_1 \times I_2 \times \cdots \times I_N}$,
$\;(k = 1, 2, \ldots, K)$, representing training data belonging to $C$ classes or categories. Each training sample $\underline{\bX}^{(k)}$
is given a label $c_k$ indicating the category (class) to which it belongs.

\begin{figure}[t]
 \centering
 \includegraphics[width=9.5 cm,height=8.2cm]{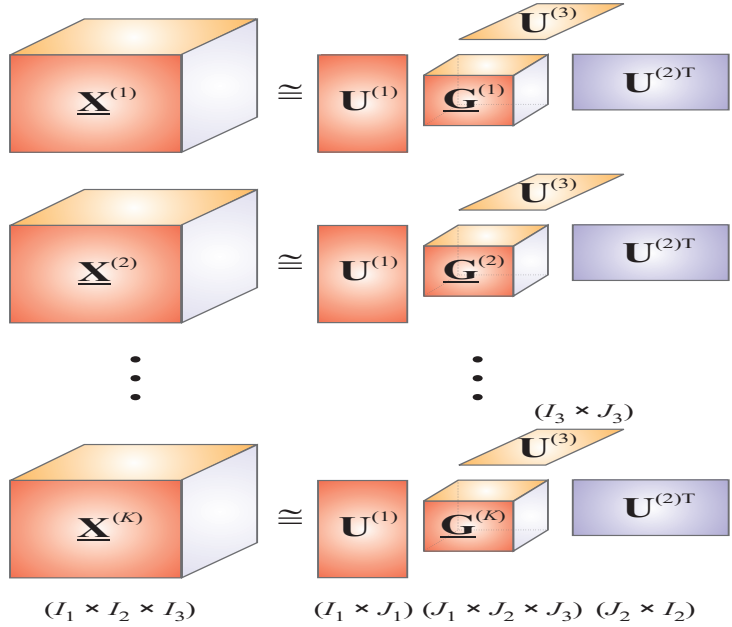}
 \caption{Illustration of feature extraction from a set of 3-way data tensors $\underline \bX^{(k)} \in \Real^{I_1 \times I_2 \times I_3}$, $(\;k=1,2,\ldots,K)$. The objective is to estimate common factor matrices (bases) $\bU^{(n)} \in \Real^{I_n \times J_n} $ $\;(n=1,2,3)$ and core tensors $\underline \bG^{(k)} \in \Real^{J_1 \times J_2  \times  J_3}$, typically with $J_n <<I_n$.}
 \label{fig:tuckersamples}
\end{figure}

In order to preform a dimensionality reduction and to extract significant features,
 we need to apply
simultaneous approximative tensor decompositions (see Fig. \ref{fig:tuckersamples})
\be
 {\underline{\bX}^{(k)}} = {\underline{\bG}^{(k)}}\, \times_1 \bU^{(1)} \, \times_2 \bU^{(2)} \cdots \, \times_N \bU^{(N)} + \underline{\bE}^{(k)},
\, 
 \label{equ_compdata}
\ee
 $(k = 1,2,\ldots,K)$, where the compressed core tensors $\underline{\bG}^{(k)} \in \Real^{J_1 \times J_2 \times \cdots \times J_N}$ representing features are of  lower dimension than the original data tensors $\underline{\bX}^{(k)}$, and we assume that the factors  (basis matrices)
$\bU^{(n)} = [\bu^{(n)}_1, \bu^{(n)}_2, \ldots, \bu^{(n)}_{J_n}] \in \Real^{I_n \times J_n}, (n= 1,2, \ldots, N)$ are common factors for all data tensors.

\begin{figure} 
 \centering
 (a)
 \includegraphics[width=13.5 cm,height=6.8cm]{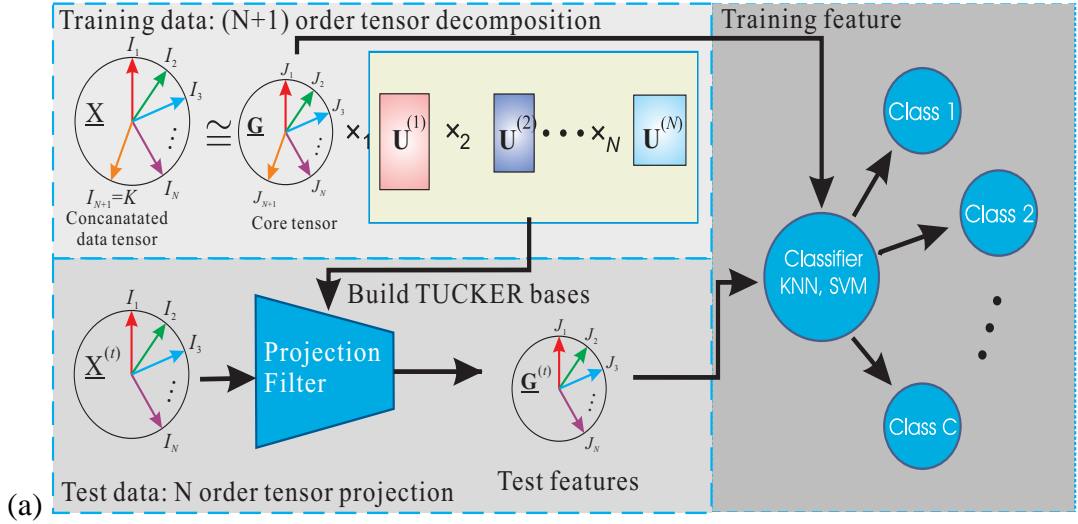}\\
 \vspace{0.7cm}
 (b)
 \includegraphics[width=13.5 cm,height=6.2cm]{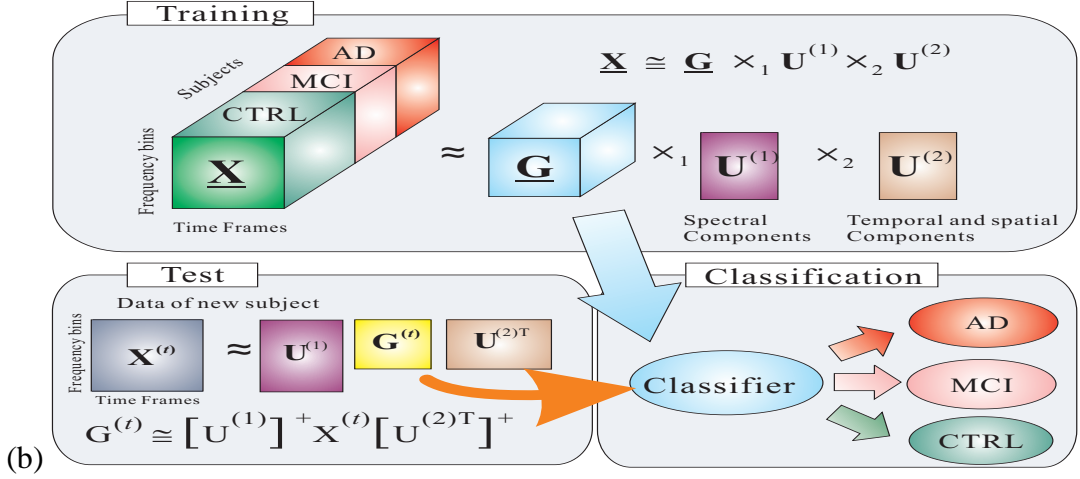}
 \caption{(a) Conceptual diagram illustrating a classification procedure
 based on the Tucker decomposition of the concatenated tensor $\underline \bX \cong {\underline \bG} \, \times_1 \bU^{(1)} \, \times_2 \bU^{(2)} \cdots \, \times_N \bU^{(N)} \in \Real^{I_1 \times I_2 \cdots I_N \times K}$ consisting of all sampling training data $\underline \bX^{(k)} \in \Real^{I_1 \times I_2 \cdots  \times I_N}$.
 Reduced features are obtained by projecting the testing data tensor $\underline \bX^{(t)}$ onto the feature subspace spanned by factors (bases) $\bU^{(n)}$. The projection filter depends on the factors  $\bU^{(n)}$. If they are orthogonal, we can apply  a simple projection  filter: $\underline \bG^{(t)} = \underline \bX^{(t)} \times_1 \bU^{(1) \; T} \times_2 \bU^{(2) \;T} \cdots  \times_N \bU^{(N) \;T}$ (otherwise, see Eq. (\ref{ProjG}). (b) Application of the general procedure for the classification of Alzheimer's disease (AD), Mild Cognitive Impairment (MCI) patients and  age matched  control subjects (CTRL).}
 \label{fig:Class}
\end{figure}

To compute the factor matrices $\bU^{(n)}$ and the core tensors $\underline \bG^{(k)}$,  we  concatenate all  training (sample) tensors into one $N+1$ order training data tensor ${\underline {\bf{X}}} = {\tt{cat}}({{\underline \bX}^{(1)}},{{\underline \bX}^{(2)},
\ldots, {{\underline \bX}^{(K)}}}, N+1) \in \Real^{I_1 \times I_2 \cdots  \times I_N \times I_{N+1}}$, with $N+1=K$ and  perform the  Tucker-$N$ decomposition \cite{Phan-Class10}:
\be
 {\underline \bX} = {\underline \bG} \, \times_1 \bU^{(1)} \, \times_2 \bU^{(2)} \cdots \, \times_N \bU^{(N)} +{\underline \bE},
\ee
where the sample tensors $\underline \bX^{(k)}$ can be extracted back from the concatenated tensor  by fixing the $(N+1)$-th index at a value $k$, i.e.,
$\underline {\bf{X}} (\mathop :\limits_1 ,\mathop :\limits_2 , \ldots ,\mathop :\limits_N ,\mathop k\limits_{N + 1} ) = \underline {\bf{X}} ^{\left( k \right)}\,$ and the
 individual features (corresponding to different classes) are extracted from the core tensor $\underline \bG \in \Real^{J_1 \times J_2 \cdots \times J_N \times J_{N+1}}$ as
$\underline {\bf{G}}^{\left( k \right)}= \underline {\bf{G}} (\mathop :\limits_1 ,\mathop :\limits_2 , \ldots ,\mathop :\limits_N ,\mathop k\limits_{N + 1} ) $, with $J_{N+1}=K$.
In other words,  the features of a specific training data $\underline{\bX}^{(k)}$ are
represented by the $k$-th row of the \mbox{mode-$(N+1)$} matricized
version of the core tensor $\underline \bG$.

The above procedure for the feature extraction has been applied for a wide class of classification problems, as illustrated in Fig \ref{fig:Class} (a) \cite{Phan-Class10}. In the first stage, we extracted features and estimated the basis matrices $\bU^{(n)}$ for the concatenated tensor of all training data $\underline \bX^{(k)}$ with labels and we built up a projection filter to extract the  features of the test data (without labels). In the next step, we extracted features of test data using estimated common  basis. The extracted features were then compared with the training features using standard classifier, e.g., LDA, KNN, SVM \cite{Phan-Class10}.
We applied the procedure described above to various feature extraction  classification problems
\cite{Phan-Class10}. For example, we applied this procedure to successfully classify three
groups of human subjects: Control age matched subjects (CTRL), Mild Cognitive Impairment (MCI)  and probable Alzheimer's disease (AD) patients (see Fig. \ref{fig:Class} (b)) on the basis of spontaneous, followup EEG data \cite{CichockiAD,Phan-Class10}. All the training EEG data were organized in the form of 3-way tensors using a time-frequency representation by applying Morlet wavelets. Each frontal slice represented preprocessed data of one subject.  We applied the standard  Tucker decomposition using HOOI and NTD algorithms to extract the basis matrices $\bU^{(n)}$  $(n=1,2)$ and reduced features represented by a core tensor $\underline \bG$. The new test data in the form of matrices $\bX^{(t)}$ were projected via the basis matrices to extract features. The extracted test features were compared with the training features and classification was performed using LDA, KNN and SVM. We obtained quite promising classification accuracies to predict AD ranging from 82\% to 96\% depending on EEG data sets.

\section{Linked Multiway BSS/ICA}

In neuroscience, we often need to perform a  so called group analysis which seeks to identify mental tasks or stimuli driven brain patterns that are common in two or more subjects in a group.

Various group ICA methods have been proposed to combine the results of multi-subjects, multi-tasks only after ICA is carried out on individual subjects and usually   no constraints are imposed
on  components, which are allowed to
differ with respect to their  spatial, spectral
maps, as well as to their temporal patterns.
However, in many scenarios some  links
need to be considered  to analyze variability and consistency of  the components
across subjects. Furthermore, some components do not need to be
 necessarily  independent, they can be instead sparse, smooth or non-negative  (e.g., for spectral components).
Moreover, it is often  necessary to  impose some constraints to be able to estimate some components, which are
 identical or maximally correlated across subjects with regard to
their spatial distributions, spectral or  temporal patterns. This leads to a new  concept and model of linked multiway BSS (or linked multiway ICA, if  statistical independence criteria are used).

\begin{figure}[t]
 \centering
 \includegraphics[width=9.5 cm,height=9.2cm]{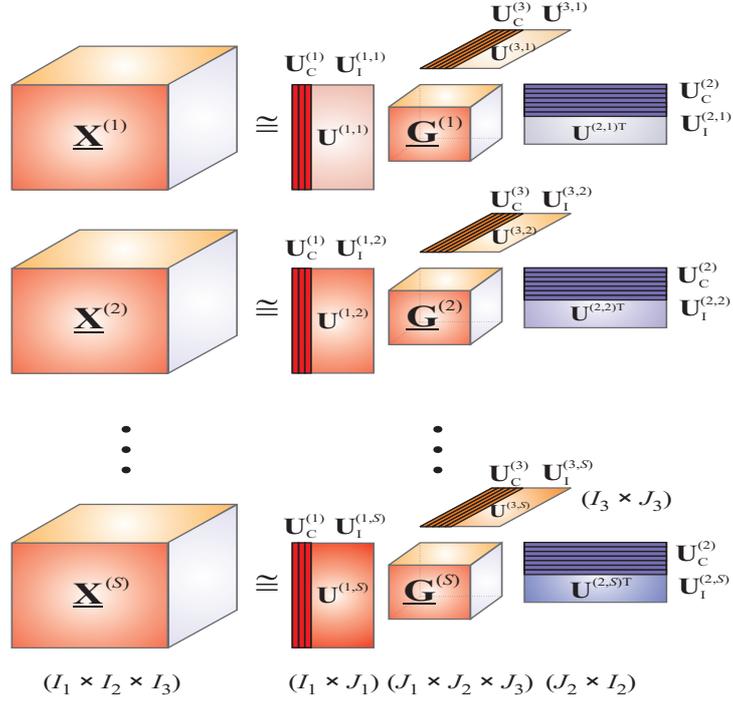}
 \caption{Conceptual model of tensors decomposition for
  linked  multiway BSS (especially, Linked Multiway ICA). Objective is to find constrained factor matrices $\bU^{(n,s)}= [\bU^{(n)}_C,\; \bU^{(n,s)}_I] \in \Real^{I_n \times J_n}$, $\;n=1,2,3)$ and core tensors $\underline \bG^{(s)} \in \Real^{J_1 \times J_2 \times J_3}$, which are partially linked or maximally correlated, i.e., they have the same common components or highly correlated components.}
 \label{fig:LMBSS}
\end{figure}

In  the linked multiway BSS, we perform approximative decompositions of a set of data tensors $ \underline \bX^{(s)} \in  \Real^{I_1 \times I_2 \cdots \times I_N}$, $\;(s=1,2,\ldots,S)$ representing multiple subjects and/or multiple tasks (see Fig.  \ref{fig:LMBSS}):
\be
   \underline \bX^{(s)} = \underline \bG^{(s)} \times_{1} \bU^{(1,s)} \times_{2} \bU^{(2,s)}   \cdots \times_{N} \bU^{(N,s)} + \underline {\bE}^{(s)},
 \label{LMBSS}
\ee
where each factor (basis matrix) $\bU^{(n,s)}=[\bU^{(n)}_C, \; \bU^{(n,s)}_I]
\in \Real^{I_n \times J_n}$
is composed of two  bases: $\bU^{(n)}_C \in \Real^{I_n \times R_n} $ (with $0 \leq R_n \leq J_n$),
  which are common bases for all subjects in the group and correspond to the same or maximally correlated  components and $\bU^{(n)}_I \in \Real^{I_n \times (J_n-R_n)}$, which correspond to stimuli/tasks independent individual characteristics.

  For example,  $\underline \bX^{(s)},\;\ (s=1,2,\ldots S)$ denotes the brain activity in space-time-frequency domains for the $s$-th subject. Based on a set of such data, we can compute   common factors and interactions between them.
 For example, linked multiway BSS approach my reveal  coupling of brain regions, possibly at  different time slots and/or different frequency bins.

To solve problems formulated this way, we can apply similar procedures to the one described in previous two Sections. If $\bU^{(n,s)}= \bU^{(n)}_C \in \Real^{I_n \times J_n}$ for a specific mode $n$, then we can concatenate all data tensors along this mode, perform tensor decomposition and  apply any standard BSS algorithm to compute desired components in each mode (e.g., to estimate independent components, we apply any standard ICA algorithm).  In the more general case, when the number of common components $R_n$ are unknown,  we perform an additional unfolding for each data tensor $\underline \bX^{(s)}$ in each mode and then perform a set of constrained matrix factorizations (by applying standard algorithms for ICA, NMF, SCA, SmoCA, PCA/SVD, etc.):
\be
\bX^{(s)}_{(n)}= \bU^{(n,s)} \; \bG_{(n)}^{(s)} \bZ^{(n,s)} + \bE^{(n)} \cong \bB_{n,s} \; \bA^T_{n,s}, 
\ee
for  $n=1,2,\ldots,N$ and  $s=1,2, \ldots,S$, where matrices $\bB_{n,s}=\bU^{(n,s)}$ represent individual linked components, while matrices $\bA^T_{n,s}=\bG_{(n)}^{(s)} \bZ^{(n,s)}$ represent basis vectors or linked mixing processes and
$\bZ^{(n,s)}=\left[\bU^{(N,s)} \otimes \cdots \otimes \bU^{(n+1,s)}
  \otimes \bU^{(n-1,s)}  \cdots \otimes \bU^{(1,s)}\right]^{T} $.

In the next stage, we need to  perform the statistical analysis and to compare the individual components extracted in each mode $n$,  for all subjects $S$, by performing clustering and ordering components to identify common  or similar (highly correlated) components (see e.g., \cite{NMF-book,Guo-GICA}).
In the last stage, we compute core tensors, which describe the functional links between components (see Eq. (\ref{ProjG})).

For the linked multiway BSS, we can also exploit constrained Block Tensor Decomposition (BTD) models for an averaging data tensor across subjects:
\be
  \underline {\bar \bX} = \sum_{s=1}^S \left(\underline \bG^{(s)} \times_{1} \bU^{(1,s)} \times_{2} \bU^{(2,s)}   \cdots \times_{N} \bU^{(N,s)}\right) + \underline {\bar \bE}.
  \label{LMBSS-BTD}
\ee
Such model can provide us additional information.

In group and linked multiway BSS and ICA, we usually seek stimuli driven ERPs (event related responses) and/or task related common
components or common bases reflecting both intra subject and inter subject features as bases, which are independent  involving individual on going
brain activities \cite{GroupNMF}.
In other words, we seek event/task-related components $\bU_C^{(n)}$, that are identical in each mode or maximally correlated across subjects and event/task independent individual bases $\bU_I^{(s,n)}$, which are independent or even as far apart as possible (anti-correlated) \cite{GroupNMF}.

Note that our Linked Multiway BSS is different from the linked ICA \cite{LinkedICA}
 and group NMF \cite{GroupNMF,Lee2007},  since we  do not limit components diversity by constraining them to  be only statistically independent or only nonnegative components and our model is based on constrained Tucker models instead of rather the quite restrictive trilinear CP model.

How to select common components  will depend on the validity of the underlying
assumptions and {\it a priori} knowledge. For instance, identical spatial distributions can
well be assumed for homogeneous subject groups, but may
be unacceptable in studies that include patients with different ages or mental diseases or
abnormalities.
Temporal components  may be the same for
stimulus- or task-evoked responses that are related to a
specific mental task or  paradigm, but these will vary for the spontaneous
fluctuations that occur in resting state experiments. In some
experiments, responses may be assumed similar or   identical within
but different across subgroups \cite{Langers-GICA,Guo-GICA}.

We  conclude that the proposed Linked BSS model
provides a framework that is  very flexible and general and it  may substantially
supplement many of the currently available  techniques for group ICA
and  feature extraction models.
Moreover, the  model can be extended to a multiway Canonical Correlation Analysis (MCCA),
 in which we impose maximal correlations among
normalized factor matrices (or subset of components) and/or core tensors.

\section{Multiway Partial Least Squares}

The Partial Least Squares (PLS) methods (originally developed in chemometrics and econometrics) are particularly suitable  for the analysis of relationships among multi-modal brain data  (e.g., EEG, MEG, ECoG (electrocorticogram), fMRI) or relationships between measures of brain activity and behavior data \cite{Abdi10,PLS-Neuro}.
The standard  PLS approaches have been recently summarized by H. Abdi \cite{Abdi10} and  their suitability to model relationships between brain activity and  behavior (experimental design) has been highlighted by A. Krishnan {\it et al.} \cite{PLS-Neuro}.

In computational neuroscience, there are two related basic  PLS methods: PLS correlation (PLSC), which analyzes correlations or associations between two or more sets of data (e.g., two modalities brain data or brain  and behavior data) and PLS regression (PLSR) methods, which attempt to predict one set of data  from another independent data  that constitutes the predictors (e.g., experimental behavior data from brain data such multichannel ECoG or  scalp EEG from ECoG,  by performing simultaneous recordings for epileptic patients).

In order to predict response variables represented by matrix $\mathbf{Y}$ from the independent variables $\mathbf{X}$, the standard PLSR techniques find a set of common orthogonal latent variables (also called latent vectors, score vectors or components) by projecting both $\mathbf{X}$ and $\mathbf{Y}$ onto a new subspace, which ensures a maximal correlation between the latent variables of $\mathbf{X}$ and $\mathbf{Y}$ \cite{Abdi10}.
In other words, the prediction is achieved by simultaneous  approximative decompositions  training data sets: $\mathbf{X} \in \Real^{I \times N} $ and  $\mathbf{Y} \in \Real^{I \times M}$ into components ($\bA=[\ba_1,\ba_2,\ldots,\ba_J] \in \Real^{I \times J} ,\bB =[\bb_1,\bb_2,\ldots,\bb_J] \in \Real^{N \times J}$ and $\bC=[\bc_1,\bc_2,\ldots,\bc_J] \in \Real^{M \times J}$) :
\be
\bX &\cong& \bA\bB^T =\sum_{j=1}^J \ba_j \bb_j^T, \\
\bY &\cong& \bA \bD \bC^T=\sum_{j=1}^J d_{jj} \ba_j \bc_j^T,
\ee
where $\bD \in \Real^{J \times J} $ is a scaling diagonal matrix; with the constraint that these components explain as much as possible of the covariance between $\bX$ and $\bY$ \cite{Abdi10,Zhao-PLS}.
The latent components $\bA$ are defined as $\bA= \bX \bW$, where $\bW =[\bw_1,\bw_2, \ldots, \bw_J] \in \Real^{N \times J}$ consists of $J$ direction vectors $\bw_j$. The basic concept of the standard PLS is to find these directions vectors $\bw_j$ from successive optimizations problems:
\be
\bw_j&=& \mbox{arg} \max_{\bw} \left\{ \mbox{corr}^2 (\bX \bw , \bY)\; \mbox{var}(\bX \bw) \right\}, \nonumber \\
s.t. \quad && \bw^T\bw = 1, \quad \bw^T \bX^T \bX \bw_k =0 \nonumber  
\ee
for $k=1,2,\ldots, j-1$.

Such decompositions and relationships between components can be used to predict values of dependent variables  for new situations, if the values of the corresponding independent variables are available.
%
%

Since the brain data  are often multidimensional and multi-modal and can be naturally represented via tensors, attempts have been made to extend PLS methods to multiway models \cite{Zhao-PLS}.
%
%
In this section, we briefly describe multiway PLS based on a constrained Tucker model (PTD):
Given an $N$th-order independent data tensor $\underline{\mathbf{X}}\in \Real^{I_{1}\times \cdots\times I_{n}\times\cdots\times I_{N} }$ and an $M$th-order dependent data tensor $\underline{\mathbf{Y}}\in \Real^{K_{1}\times \cdots \times K_{m} \times\cdots\times K_{M}}$, with the same size in at least   one  mode (typically, the first mode)\footnote{ The first mode,  usually, corresponds to the samples mode or the time mode. However,  in some applications, for example,  for simultaneous recordings of   EEG and ECoG data two  modes (ways) may  have the same size and they represent time and frequency (temporal and spectral modes, respectively).} $I_1 = K_1$.

\begin{figure}[t]
\centering
\includegraphics[width=11.1 cm,height=8.2cm]{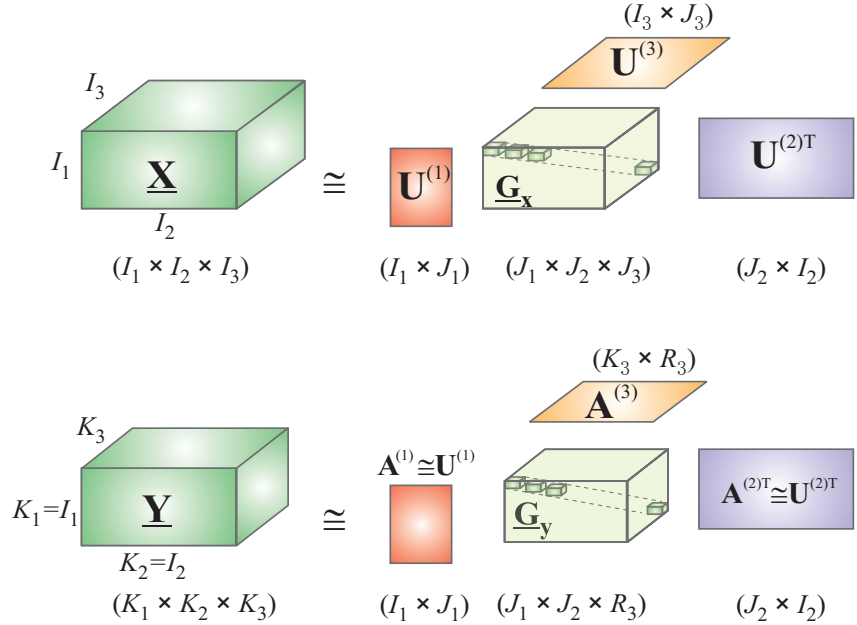}
\caption{Illustration of a generalized  Tucker-based Multiway PLS model for 3-way tensors in the case when two factor matrices are common (see Eqs. (\ref{HOPLS1})-(\ref{HOPLS2})). Our objective is to perform  approximative decompositions  of independent tensor $\underline{\mathbf{X}} \in \Real^{I_1 \times I_2 \cdots \times I_N}$ and dependent tensor $\underline{\mathbf{Y}} \in \Real^{K_1 \times K_2 \cdots \times K_M} $, with $I_1=K_1$ and $I_2=K_2$, by imposing  additional constraints, especially that the two factor matrices  are the same or highly correlated: $\bA^{(1)} \cong \bU^{(1)}$ and $\bA^{(2)} \cong \bU^{(2)}$, while the other factors are uncorrelated (or  anti-correlated) or statistically independent. We assumed also that the core tensors are sparse and block diagonal.}
\label{fig:HOPLS}
\end{figure}

Our objective is to preform simultaneous constrained  Tucker decomposition, with  at least one (or two) common or maximally correlated factor(s) (see Fig. \ref{fig:HOPLS}):
\be
\label{HOPLS1}
   \underline \bX &=& \underline \bG_{\bx} \times_{1} \bU^{(1)} \times_{2} \bU^{(2)}\times_{3}  \cdots \times_{N} \bU^{(N)} + \underline {\bE}_{\bx}, \\
      \underline \bY &=& \underline{\bG}_{\by} \times_{1} \bA^{(1)} \times_{2} \bA^{(2)} \times_{3} \cdots \times_M \bA^{(M)}+\underline {\bE}_{\by},
      \label{HOPLS2}
\ee
where additional constraints are imposed:  $\bU^{(1)}\cong\bA^{(1)} \in \Real^{I_1 \times J_1}$  (and $\bU^{(2)} \cong \bA^{(2)} \in \Real^{I_2 \times J_2}$ of $I_2=K_2$), while other factor matrices are essentially different (e.g., orthogonal or mutually independent).
 %
 Note that the core tensors $\underline {\bG}_{\bx}$ and $\underline {\bG}_{\by}$ have  special block-diagonal structures (see Fig. \ref{fig:HOPLS}) that indicate sparseness. New algorithms for constrained  Tucker based  multiway PLS models  have been developed in \cite{Zhao-PLS}.

Such models allow for different types of structures on$\underline \bX$ and $\underline \bY$ and provides a general framework for solving multiway regression problems that explore complex relationships between multidimensional dependent and independent variables.
%
%
For example, tensor decompositions  can be  applied in emerging neuroimaging genetics  studies  to investigate links between biological parameters measured with brain imaging and genetic variability \cite{Imaging-genetics}.
%

%
\label{sec:HOPLS}

\section{Conclusions}

Tensor decompositions is a fascinating emerging field of research, with
 many applications in multilinear blind source separation, linked multiway BSS/ICA, feature extraction, classification and prediction. In this
review/tutorial paper, we have briefly discussed several new and promising models for  decompositions of linked or correlated sets of tensors, by imposing various constraints on factors, components or latent variables depending on the specific applications.
We have developed several  promising and efficient algorithms for such models,
 which can be found in our recent book, publications and reports.

\bibliographystyle{ieicetr}

\end{document}